\documentclass{article}

\usepackage{arxiv}

\usepackage{amsmath, amssymb, amsthm}
\usepackage{graphicx}    % include graphics
\usepackage{subcaption}  % subfigures (recommended)
\usepackage{float}       % [H] placement if needed

\usepackage{geometry}
\usepackage{graphicx}

\usepackage{setspace}
\usepackage{microtype}

\usepackage[colorlinks=true,
            linkcolor=blue,
            citecolor=blue,
            urlcolor=blue]{hyperref}

\title{Data-Driven Discovery of Sign-Indefinite Artificial Viscosity for Linear Convection\\
\large A Space--Time Reconvolution Perspective}

\author{ \href{}{\hspace{1mm}Arun Govind Neelan} \\
	SimuNetics\\
	Kaniyakumari\\
	India-629173 \\
	\texttt{arunneelaniist@gmail.com}
}

\renewcommand{\shorttitle}{Data-Driven Sign-Indefinite Viscosity}

\hypersetup{
pdftitle={Data-Driven Discovery of Sign-Indefinite Artificial Viscosity},
pdfsubject={Computational Physics, Machine Learning},
pdfauthor={Arun Govind Neelan},
pdfkeywords={Artificial viscosity, Data-driven closure, CFD, Scientific Machine Learning},
}

\begin{document}
\maketitle

\begin{abstract}
Artificial viscosity is traditionally interpreted as a positive, spatially acting regularization introduced to stabilize numerical discretizations of hyperbolic conservation laws.
In this work, we report a data-driven discovery that motivates a reinterpretation of this classical view.
We consider the linear convection equation discretized using an unstable FTCS scheme augmented with a learnable artificial viscosity.
Using automatic differentiation and gradient-based optimization, the viscosity field is inferred by minimizing the error with respect to the exact solution, without imposing any sign constraints.
The optimized viscosity consistently becomes \emph{locally negative} near extrema, while the numerical solution remains stable and nearly exact.
This behavior is not readily explained within classical modified equation analysis and Lax--Wendroff-type arguments, which predict a strictly positive effective viscosity.
To resolve this apparent contradiction, we reinterpret artificial viscosity as a space–time closure that compensates unresolved truncation errors while enforcing entropy stability through global dissipation balance rather than pointwise positivity.
Within this framework, the Lax--Wendroff scheme corresponds to a degenerate projection in which temporal truncation errors are eliminated and reintroduced as spatial diffusion.
We show that entropy stability constrains the integrated dissipation budget rather than the pointwise sign of spatial viscosity.
As a result, locally negative viscosity naturally emerges as a numerical reconvolution operator that compensates for dispersive truncation errors.
Negative viscosity is therefore not an unphysical diffusion process, but a scheme- and grid-dependent correction mechanism.
\end{abstract}


\section{Introduction}

Numerical discretizations of hyperbolic conservation laws inevitably introduce truncation errors that manifest as artificial dissipation and dispersion.
To control these effects, artificial viscosity has long been employed as a stabilization mechanism \cite{vonNeumann1950}, typically interpreted as a positive, spatially acting regularization motivated by entropy arguments and monotonicity requirements \cite{Tadmor1987}.
In classical numerical analysis, viscosity is most often viewed as a spatial correction added to counteract instabilities arising from discrete advection operators \cite{Jameson1981}.
Within this perspective, time integration errors are either neglected or eliminated via Taylor expansion and subsequently reintroduced as spatial diffusion, as in Lax--Wendroff \cite{Lax1960} and Taylor--Galerkin schemes \cite{Donea1984}.
While this approach has proven highly successful, leading to Total Variation Diminishing (TVD) \cite{Harten1983} and Essentially Non-Oscillatory (ENO) schemes \cite{Shu1988}, it implicitly assumes how truncation errors should be redistributed between space and time.

Recent advances in data-driven and optimization-based PDE solvers provide an opportunity to revisit these assumptions.
Techniques ranging from neural networks to differentiable programming have been used to improve dispersion relation analysis \cite{NeelanDispersion2025} and discover closure terms that are consistent with data \cite{BarSinai2019}. Sengupta et al.~\cite{sengupta2023quantifying} rigorously analyze the explicit central-difference Lax–Wendroff method using global spectral analysis to determine optimal numerical parameter ranges that achieve prescribed accuracy for high-fidelity simulations of 2D convection–diffusion problems.
The emergence of Physics-Informed Neural Networks (PINNs) \cite{Raissi2019, neelan2024physics} and differentiable programming for fluid dynamics \cite{Holl2020, Kochkov2021} has enabled the discovery of closure terms that may be complex in form.
Rather than prescribing the form and sign of artificial viscosity \emph{a priori} \cite{Duraisamy2019}, one may infer stabilizing closures directly by minimizing solution error with respect to reference data \cite{Sirignano2018}.

In this work, we apply gradient-based optimization with automatic differentiation \cite{Baydin2018} to the linear convection equation discretized using an unstable FTCS scheme augmented with learnable artificial viscosity.
Without imposing sign or monotonicity constraints, the optimized viscosity field consistently develops locally negative values near smooth extrema, while the numerical solution remains stable and nearly exact.
This behavior is not readily explained within classical modified equation analysis \cite{Warming1974}, which predicts a strictly positive effective viscosity.

The objective of this paper is to show that this apparent discrepancy arises from a purely spatial interpretation of artificial viscosity.
We propose instead to view viscosity as a \emph{space--time closure} for unresolved truncation errors introduced by discretization.
Within this framework, the Lax--Wendroff scheme emerges as a degenerate space--time projection in which temporal truncation errors are eliminated and reintroduced as spatial diffusion.
Once dissipation is allowed to redistribute across space and time, locally negative spatial viscosity becomes admissible and, in some cases, necessary.

We further show that entropy stability constrains the integrated dissipation budget rather than the pointwise sign of spatial viscosity \cite{Guermond2011}.
As a consequence, negative viscosity is naturally interpreted as a numerical reconvolution operator that compensates dispersive truncation errors, rather than as an unphysical diffusion process.
This perspective unifies data-driven closures, Taylor--Galerkin stabilization, spectral viscosity \cite{Tadmor1989}, and compact flux corrections \cite{Lele1992} under a single space--time reconvolution framework.

\section{Data-Driven Discovery of Sign-Indefinite Viscosity}

In this section, we describe the data-driven optimization procedure that led to the empirical observation motivating this work.
The objective is not to propose a new numerical scheme, but to infer a stabilizing closure for a fixed, consistent discretization by minimizing solution error with respect to reference data.

\subsection{Baseline Discretization}

We consider the linear convection equation
\begin{equation}
u_t + a u_x = 0 ,
\label{eq:linconv_dd}
\end{equation}
posed on a periodic domain.
The equation is discretized using a forward-time, centered-space (FTCS) scheme,
\begin{equation}
\frac{u_i^{n+1} - u_i^n}{\Delta t}
+
a \frac{u_{i+1}^n - u_{i-1}^n}{2\Delta x}
= 0 .
\label{eq:ftcs}
\end{equation}
Although consistent, this scheme is linearly unstable and therefore requires stabilization.
To regularize the discretization, we augment \eqref{eq:ftcs} with an artificial viscosity term written in conservative form,
\begin{equation}
\frac{u_i^{n+1} - u_i^n}{\Delta t}
+
a \frac{u_{i+1}^n - u_{i-1}^n}{2\Delta x}
=
\frac{1}{\Delta x}
\left[
\mu_{i+1/2} (u_{i+1}^n - u_i^n)
-
\mu_{i-1/2} (u_i^n - u_{i-1}^n)
\right],
\label{eq:ftcs_visc}
\end{equation}
where $\mu_{i+1/2}$ denotes a cell-face viscosity coefficient.
No assumptions are made regarding the sign or smoothness of $\mu$.

\subsection{Optimization Problem}

The viscosity field $\mu$ is inferred by minimizing the discrepancy between the numerical solution and the exact solution of \eqref{eq:linconv_dd}.
Let $u_{\mathrm{num}}(x,t;\mu)$ denote the solution produced by \eqref{eq:ftcs_visc}, and let $u_{\mathrm{exact}}(x,t)$ denote the analytical solution.
We define the objective function
\begin{equation}
\mathcal{J}(\mu)
=
\frac{1}{2}
\int_0^T
\int_\Omega
\left(
u_{\mathrm{num}}(x,t;\mu)
-
u_{\mathrm{exact}}(x,t)
\right)^2
\, dx \, dt .
\label{eq:loss}
\end{equation}

The optimization problem is then
\begin{equation}
\mu^\ast
=
\arg\min_{\mu}
\mathcal{J}(\mu),
\label{eq:opt}
\end{equation}
subject only to the discrete time-marching scheme \eqref{eq:ftcs_visc}.
No positivity, monotonicity, or entropy constraints are imposed.
Gradients of $\mathcal{J}$ with respect to $\mu$ are computed using automatic differentiation through the fully discrete time integration.
All results reported here are obtained using gradient-based optimization.

\subsection{Observed Structure of the Optimized Viscosity}

Across a range of grid resolutions, time step sizes, and smooth initial conditions, the optimized viscosity field $\mu^\ast$ exhibits a consistent qualitative structure.
In regions where the solution is locally monotone, $\mu^\ast$ remains small and positive.
In contrast, near under-resolved transition regions where the solution departs from piecewise smoothness, the optimized viscosity develops locally negative values. Despite the presence of negative viscosity, the numerical solution remains stable, exhibits no spurious oscillations, and closely matches the exact solution in both phase and amplitude.
The total entropy of the system is observed to remain non-increasing in time.
This behavior persists under mild regularization of $\mu$ and is robust to changes in initial conditions.
However, constraining $\mu$ to remain non-negative leads to systematic amplitude decay near smooth extrema, consistent with classical diffusive behavior.

\subsection{Methodological Remarks and Scope}

\subsubsection{Discrete versus Continuous Adjoint}

The optimization problem considered in this work seeks to minimize the error of a fully discrete numerical solution.
For this reason, gradients are computed with respect to the discrete time-marching scheme rather than its continuous counterpart \cite{Giles2000}.
While continuous adjoint formulations are well suited for sensitivity analysis of continuous partial differential equations \cite{Jameson1988}, their use in the present context would require an additional level of approximation associated with discretizing the adjoint equations.
In contrast, the discrete adjoint obtained via automatic differentiation provides exact gradients of the discrete objective function with respect to the discrete viscosity coefficients.

\subsubsection{Why Not Neural Networks?}

Neural networks are frequently employed in data-driven PDE solvers to represent unknown closures or subgrid-scale models \cite{Beck2019, Pathak2022}.
However, the objective of the present work is not to learn a black-box surrogate, but to identify the minimal numerical correction required to stabilize a fixed, consistent discretization.
By parameterizing artificial viscosity directly at the discrete level, the present approach avoids introducing additional representation or training complexity and allows the optimizer to focus solely on compensating truncation errors associated with the underlying space--time discretization.

\section{Space--Time Closure Interpretation}

The data-driven experiments presented in Section~2 indicate that artificial viscosity inferred through error minimization does not behave as a purely dissipative spatial operator.
In particular, the emergence of locally negative viscosity near smooth extrema suggests that classical interpretations based solely on spatial stabilization are incomplete.
In this section, we propose an alternative interpretation in which artificial viscosity is viewed as a space–time closure for discretization-induced truncation errors, while enforcing numerical stability through a globally constrained dissipation mechanism.
This perspective provides a unifying framework for understanding the observed optimization results and clarifies their relationship to classical schemes.

\subsection{From Spatial Regularization to Space--Time Closure}

Consider the linear convection equation
\begin{equation}
u_t + a u_x = 0 .
\label{eq:linconv_st}
\end{equation}
In classical numerical analysis, stabilization is typically introduced through a purely spatial diffusion term,
\begin{equation}
u_t + a u_x
=
\partial_x \left( \mu u_x \right),
\label{eq:spatial_visc}
\end{equation}
where $\mu \ge 0$ is interpreted as an artificial viscosity coefficient.
This formulation implicitly assumes that truncation errors arise primarily from the spatial discretization.
However, fully discrete schemes introduce errors through both space and time approximations \cite{Deconinck2013}.
When temporal truncation errors are eliminated analytically or implicitly absorbed into spatial operators, their influence is no longer explicit in the governing equations.
To account for this, we consider a generalized space--time regularization of the form
\begin{equation}
u_t + a u_x
=
\partial_x \left( \mu_s u_x \right)
+
\partial_t \left( \mu_t u_t \right),
\label{eq:space_time_closure}
\end{equation}
where $\mu_s(x,t)$ and $\mu_t(x,t)$ represent spatial and temporal closure coefficients, respectively.
This formulation does not introduce new physics.
Rather, it provides a bookkeeping mechanism for redistributing unresolved truncation errors between space and time in a controlled manner.

\subsection{Interpretation of Learned Viscosity}

Within the space--time closure framework \eqref{eq:space_time_closure}, the viscosity field inferred in Section~2 should not be interpreted as a physical diffusion coefficient.
Instead, it represents the projection of unresolved space--time truncation errors onto a reduced set of admissible operators.
In particular, when the temporal closure term $\partial_t(\mu_t u_t)$ is eliminated or absorbed into the spatial operator through Taylor expansion, its leading-order effect appears as an effective spatial diffusion term.
Depending on the local structure of the solution and the nature of the truncation error being compensated, the resulting spatial operator need not be pointwise positive.

\section{Time--to--Space Projection and the Lax--Wendroff Scheme}

In this section, we revisit the classical Lax--Wendroff scheme through the lens of the space--time closure framework introduced in Section~3.
Our objective is not to re-derive the scheme, but to reinterpret its stabilizing mechanism as a projection of temporal truncation errors onto spatial operators.
While recent work has quantified the parameter ranges required to maintain high fidelity and prescribed accuracy using this method, standard derivations typically assume a fixed relationship between temporal and spatial errors.

\subsection{Second-Order Time Expansion}

Consider the linear convection equation
\begin{equation}
u_t + a u_x = 0 .
\label{eq:linconv_lw}
\end{equation}
A second-order Taylor expansion in time yields
\begin{equation}
u^{n+1}
=
u^n
+
\Delta t\, u_t^n
+
\frac{\Delta t^2}{2} u_{tt}^n
+
\mathcal{O}(\Delta t^3).
\label{eq:taylor_time}
\end{equation}
Using \eqref{eq:linconv_lw}, the temporal derivatives may be expressed as
\begin{equation}
u_t = -a u_x,
\qquad
u_{tt} = a^2 u_{xx}.
\label{eq:time_derivatives}
\end{equation}
Substitution into \eqref{eq:taylor_time} yields the semi-discrete update
\begin{equation}
u^{n+1}
=
u^n
-
a \Delta t\, u_x^n
+
\frac{a^2 \Delta t^2}{2} u_{xx}^n
+
\mathcal{O}(\Delta t^3).
\label{eq:lw_update}
\end{equation}

\subsection{Modified Equation Interpretation}

Dividing \eqref{eq:lw_update} by $\Delta t$ and rearranging terms leads to the modified equation
\begin{equation}
u_t + a u_x
=
\frac{a^2 \Delta t}{2} u_{xx}
+
\mathcal{O}(\Delta t^2).
\label{eq:lw_modified}
\end{equation}
Equation \eqref{eq:lw_modified} is commonly interpreted as the addition of a positive artificial viscosity with coefficient
\begin{equation}
\mu_{\mathrm{LW}} = \frac{a^2 \Delta t}{2}.
\label{eq:lw_viscosity}
\end{equation}

Within the space--time closure framework \eqref{eq:space_time_closure}, however, \eqref{eq:lw_modified} admits an alternative interpretation.
Specifically, the Lax--Wendroff scheme corresponds to the degenerate choice
\begin{equation}
\mu_s = 0,
\qquad
\mu_t = \frac{\Delta t}{2},
\label{eq:lw_closure}
\end{equation}
for which the temporal closure term
$\partial_t(\mu_t u_t)$ is eliminated via Taylor expansion and reintroduced as a spatial second-derivative operator \cite{Hirsch2007}.
Thus, the effective viscosity \eqref{eq:lw_viscosity} does not arise from an explicit spatial stabilization mechanism, but from a time-to-space projection of temporal truncation error.

\section{Entropy Stability in Space--Time}

The reinterpretation of artificial viscosity as a space--time closure raises natural questions regarding stability and entropy consistency.
In particular, the appearance of locally negative spatial viscosity may seem at odds with classical entropy arguments, which are often invoked to justify pointwise positivity of diffusive operators.

\subsection{Entropy Balance for Space--Time Regularization}

Consider the space--time regularized conservation law
\begin{equation}
u_t + a u_x
=
\partial_x(\mu_s u_x)
+
\partial_t(\mu_t u_t),
\label{eq:st_entropy}
\end{equation}
where $\mu_s(x,t)$ and $\mu_t(x,t)$ denote spatial and temporal closure coefficients.
For the quadratic entropy $\eta(u) = \frac{1}{2} u^2$,
multiplication of \eqref{eq:st_entropy} by $u$ and integrating over the space--time domain (assuming periodic or compact support) yields:
\begin{equation}
\int_0^T \int_\Omega
\left(
\mu_s (u_x)^2
+
\mu_t (u_t)^2
\right)
\, dx \, dt
=
-
\int_0^T \int_\Omega
\left(
(\partial_x \mu_s) u u_x
+
(\partial_t \mu_t) u u_t
\right)
\, dx \, dt .
\label{eq:entropy_integral}
\end{equation}

Equation \eqref{eq:entropy_integral} demonstrates that entropy stability is governed by a balance between dissipative and compensatory terms in space and time \cite{Tadmor2003}.
In particular, pointwise positivity of $\mu_s$ is sufficient but not necessary for global entropy stability when temporal closure is present.

\subsection{Interpretation of Negative Spatial Viscosity}

Within the space--time closure framework, negative spatial viscosity does not necessarily imply negative entropy production.
Rather, it reflects a redistribution of dissipation between spatial and temporal channels.
Near smooth extrema, where $u_x = 0$, the leading truncation error of high-order discretizations is often dispersive.
In such regions, a locally negative spatial viscosity can act as a reconvolution operator that compensates dispersive errors, while temporal closure terms or higher-order operators maintain global stability.

\section{Implicit Emergence of Sign-Indefinite Viscosity in Classical Schemes}

The appearance of locally negative viscosity in the data-driven closures discussed in Sections~2 and~5 may initially seem unconventional.
However, sign-indefinite correction operators have long appeared implicitly in classical numerical schemes designed to balance dissipation and dispersion.

\subsection{Connection to Tadmor-Type Entropy-Stable Central Schemes}

The emergence of sign-indefinite spatial viscosity in the present data-driven framework
is closely related to dissipation mechanisms introduced in the entropy-stable central
schemes developed by Tadmore et al.~\cite{Tadmor2003}.
These schemes were designed to stabilize non-dissipative central discretizations by
adding the minimum amount of numerical dissipation required to satisfy a discrete
entropy inequality.

For a scalar conservation law
\begin{equation}
\frac{\partial u}{\partial t} + \frac{\partial f(u)}{\partial x} = 0,
\end{equation}
Tadmor’s entropy-stable central flux can be written in semi-discrete form as
\begin{equation}
\frac{\mathrm{d} u_i}{\mathrm{d} t}
=
-\frac{1}{\Delta x}
\left(
f_{i+1/2}^{\mathrm{EC}} - f_{i-1/2}^{\mathrm{EC}}
\right)
+
\frac{1}{\Delta x}
\left(
D_{i+1/2}(u_{i+1}-u_i)
-
D_{i-1/2}(u_i-u_{i-1})
\right),
\label{eq:tadmor_central}
\end{equation}
where $f^{\mathrm{EC}}$ denotes an entropy-conservative flux and $D_{i+1/2}$ is a
solution-dependent dissipation operator.
Importantly, the dissipation matrix $D_{i+1/2}$ is not equivalent to a constant positive
viscosity. Instead, it vanishes in smooth regions and activates only where required by
entropy stability. A modified-equation analysis reveals that the resulting truncation
error does not correspond to a uniformly positive second-derivative operator, but
includes higher-order derivative contributions that may locally counteract diffusion.
Thus, when expressed in physical space, the effective low-order dissipation may exhibit
sign-indefinite behavior while the full scheme remains entropy stable.

\subsection{Entropy-Conservative Core with Entropy-Stable Correction}

A related and more general formulation introduced by Fjordholm, Siddhartha Mishra and Tadmor,
the numerical flux as~\cite{Fjordholm2012}
\begin{equation}
f_{i+1/2}
=
f_{i+1/2}^{\mathrm{EC}}
-
f_{i+1/2}^{\mathrm{ES}},
\label{eq:ec_es_split}
\end{equation}
where $f^{\mathrm{EC}}$ is an entropy-conservative flux and $f^{\mathrm{ES}}$ is a
dissipative correction chosen to enforce entropy stability.

The entropy-stable correction typically takes the form
\begin{equation}
f_{i+1/2}^{\mathrm{ES}}
=
\frac{1}{2}
\mathbf{R}_{i+1/2}
\boldsymbol{\Lambda}_{i+1/2}
\mathbf{R}_{i+1/2}^{-1}
\left( u_{i+1} - u_i \right),
\label{eq:entropy_diss}
\end{equation}
where $\boldsymbol{\Lambda}$ contains characteristic wave speeds and $\mathbf{R}$ is the
eigenvector matrix. While this operator guarantees non-negative entropy production, its representation in physical space does not correspond to a simple positive viscosity
term. When expanded in a modified-equation sense, the leading second-derivative
term may locally change sign due to the interaction with higher-order terms.
This interpretation closely aligns with the learned closures presented in this work,
where a nominally unstable central discretization is stabilized by a data-driven,
space--time-dependent correction that redistributes dissipation without enforcing
uniform positivity.

\subsection{Relation to OpenFOAM-Type Central Discretizations}

It is worth noting that this philosophy is implicitly adopted in many finite-volume
solvers used in practice, including the open-source CFD framework.
Several OpenFOAM solvers employ central or central-upwind fluxes augmented with
solution-dependent scalar dissipation or flux limiters rather than explicit Riemann
solvers. Although typically described in terms of boundedness or numerical diffusion,
these stabilization mechanisms are conceptually equivalent to Tadmor-type entropy-stable
corrections.

From this perspective, the sign-indefinite artificial viscosity learned in the present
work may be viewed as an explicit, data-driven realization of dissipation mechanisms
that are already embedded implicitly in widely used central finite-volume schemes.

\subsection{Spectral Viscosity Methods}

Spectral viscosity methods \cite{Tadmor1989} were introduced to stabilize high-order spectral approximations while preserving accuracy in smooth regions.
The central idea is to damp unresolved high-frequency modes without introducing excessive dissipation at resolved scales.
Consider the linear convection equation augmented with a spectral viscosity operator.
In Fourier space, the dissipation introduced is strictly non-negative and scale-selective.
However, when the operator is expressed in physical space and formally expanded in lower-order derivatives, it admits a representation of the form
\begin{equation}
\varepsilon_N (-1)^{s+1} \partial_x^{2s} u
\sim
\mu_2 u_{xx}
+
\mu_4 u_{xxxx}
+
\mu_6 u_{xxxxxx}
+
\cdots ,
\label{eq:spectral_reconv}
\end{equation}
where the coefficients $\mu_{2k}$ alternate in sign.
As a result, the effective second-derivative contribution may be locally negative, even though the full operator remains globally dissipative.

\subsection{Compact and Dispersive Flux Corrections}

Compact and corrected central schemes introduce higher-order flux corrections to improve accuracy near smooth extrema \cite{Lele1992}.
These corrections are often designed to counteract excessive numerical diffusion introduced by first-order stabilizing terms.
A modified equation analysis typically yields a term of the form $\gamma \Delta x^2 u_{xxxx}$.
Near smooth extrema, this fourth-derivative term counteracts the leading diffusive contribution, effectively producing a locally negative second-derivative correction.
Although this behavior is rarely described in terms of negative viscosity, it reflects the same reconvolution mechanism identified in the data-driven closures of Section~2.

\section{Results and Discussion}

We consider the one-dimensional linear advection equation
\begin{equation}
\frac{\partial u}{\partial t} + c \frac{\partial u}{\partial x} = 0,
\label{eq:advection}
\end{equation}
posed on a periodic domain $x \in [0,1)$ with constant advection speed $c=1$. The domain is discretized using $N=100$ uniform grid points with $\Delta x=0.01$,
and the solution is advanced in time using $\Delta t=10^{-3}$,
corresponding to a CFL number of $0.1$, which is very diffusive for the first-order upwind schemes.

\subsection{Initial condition and exact solution}

The initial condition is a discontinuous hat function given by
\begin{equation}
u(x,0) =
\begin{cases}
1, & 0.4 < x < 0.6, \\
0, & \text{otherwise},
\end{cases}
\label{eq:ic}
\end{equation}
with periodic boundary conditions.  
The exact solution is a pure translation of the initial condition,
\begin{equation}
u(x,t) = u(x-ct).
\label{eq:exact}
\end{equation}

\subsection{Numerical scheme with learned artificial viscosity}

% A forward-time central-space (FTCS) discretization is employed,
% \begin{equation}
% u_i^{n+1} = u_i^n
% - \frac{c \Delta t}{2\Delta x} \left( u_{i+1}^n - u_{i-1}^n \right)
% + \mu_i^n \frac{\Delta t}{\Delta x^2}
% \left( u_{i+1}^n - 2u_i^n + u_{i-1}^n \right),
% \label{eq:ftcs_mu}
% \end{equation}
% where $\mu_i^n$ denotes a space--time-dependent artificial viscosity.

A forward-time central-space (FTCS) discretization is employed.  
Introducing left- and right-biased artificial viscosities 
$\mu_i^{-,n}$ and $\mu_i^{+,n}$ at the cell interfaces, the scheme
can be written in conservative form as
\begin{equation}
u_i^{n+1}
=
u_i^n
-
\frac{\Delta t}{\Delta x}
\left(
F_{i+1/2}^n - F_{i-1/2}^n
\right),
\label{eq:ftcs_conservative}
\end{equation}
where the numerical flux is defined by
\begin{equation}
F_{i+1/2}^n
=
c\,\frac{u_{i+1}^n + u_i^n}{2}
-
\frac{\mu_{i+1/2}^n}{\Delta x}
\left(
u_{i+1}^n - u_i^n
\right),
\label{eq:ftcs_flux}
\end{equation}
with
\begin{equation}
\mu_{i+1/2}^n = \mu_i^{+,n},
\qquad
\mu_{i-1/2}^n = \mu_i^{-,n}.
\end{equation}

Expanding the flux difference yields
\begin{equation}
\begin{aligned}
u_i^{n+1}
=
u_i^n
&-
\frac{c \Delta t}{2\Delta x}
\left(
u_{i+1}^n - u_{i-1}^n
\right)
\\[4pt]
&+
\frac{\Delta t}{\Delta x^2}
\Big[
\mu_i^{+,n}\left(u_{i+1}^n - u_i^n\right)
-
\mu_i^{-,n}\left(u_i^n - u_{i-1}^n\right)
\Big].
\end{aligned}
\label{eq:ftcs_mu_pm}
\end{equation}

When $\mu_i^{+,n} = \mu_i^{-,n} = \mu_i^n$, the scheme reduces to the
standard FTCS discretization with space--time-dependent artificial
viscosity,
\begin{equation}
\label{eq:non_con}
u_i^{n+1} = u_i^n
- \frac{c \Delta t}{2\Delta x}
\left( u_{i+1}^n - u_{i-1}^n \right)
+ \mu_i^n \frac{\Delta t}{\Delta x^2}
\left( u_{i+1}^n - 2u_i^n + u_{i-1}^n \right).
\end{equation}

The conservative flux formulation
\eqref{eq:ftcs_conservative}--\eqref{eq:ftcs_flux}
is closely related in structure to Tadmor-type central schemes \eqref{eq:tadmor_central}, which
are known to provide robust shock-capturing capabilities in high-speed
compressible flows, although such applications are not considered in the present work.
In contrast, the formulation in \eqref{eq:non_con} is written in
non-conservative form, while this may not pose significant difficulties for low-speed flows, it can lead to incorrect shock speeds if the artificial viscosity is not properly tuned.

In particular, the numerical flux consists of a central (entropy-neutral)
advective flux augmented by a symmetric interface dissipation term.
The interface viscosities $\mu_{i\pm1/2}^n$ play a role analogous to
The numerical viscosity in Tadmor’s entropy-stable fluxes~\cite{Tadmor2003} controls the amount of dissipation added to stabilize the central discretization.
Unlike classical Tadmor schemes, however, the viscosities
$\mu_i^{+,n}$ and $\mu_i^{-,n}$ are allowed to vary independently in
space and time. As a result, the present formulation should be interpreted as a
Tadmor-inspired space--time closure.
Here $\mu_i^n$
is treated as a trainable variable and optimized at each time step by minimizing the
instantaneous error
\begin{equation}
\mathcal{L}^n = \frac{1}{N} \sum_{i=1}^N
\left( u_i^{n+1} - u_i^{\text{exact}}(t^{n+1}) \right)^2.
\label{eq:loss}
\end{equation}
Gradient descent is applied locally in space to update $\mu_i^n$, subject to
stability constraints,
\begin{equation}
\mu_{\min} \le \mu_i^n \le \mu_{\max}.
\end{equation}

\subsection{Learned space--time structure of artificial viscosity}

\begin{figure}[H]
\centering

\begin{subfigure}{0.48\textwidth}
    \centering
    \includegraphics[height=4.5cm]{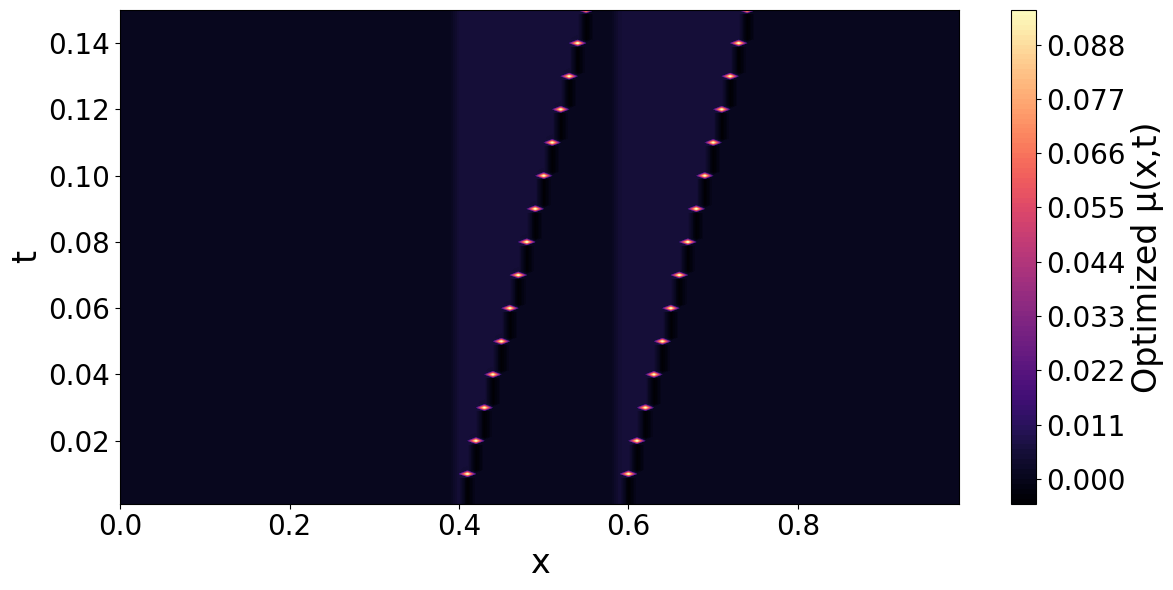}
    \caption{Space--time evolution of learned viscosity $\mu(x,t)$.}
    \label{fig:mu_xt}
\end{subfigure}
\hfill
\begin{subfigure}{0.48\textwidth}
    \centering
    \includegraphics[height=4.5cm]{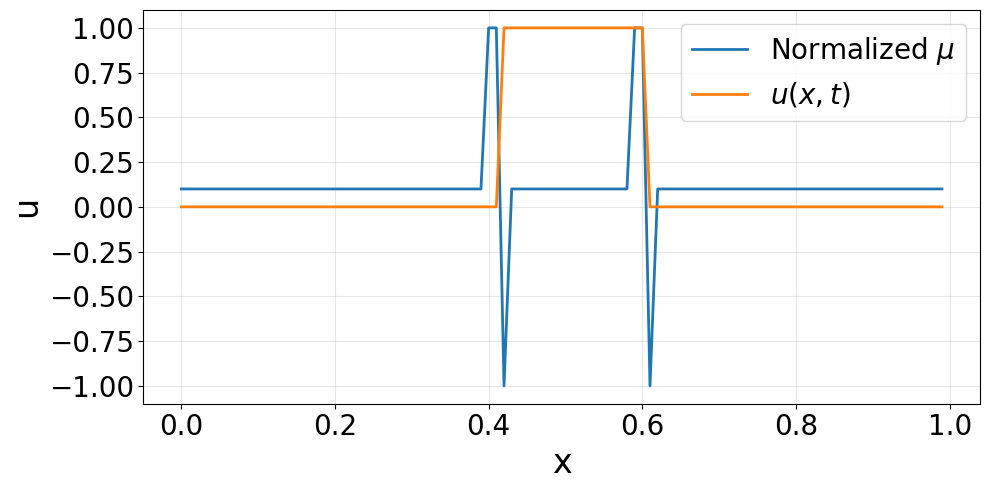}
    \caption{Learned normalized $\mu$  at T = 0.15}
    \label{fig:mu_at_T}
\end{subfigure}

\caption{Comparison of learned artificial viscosity }
\label{fig:mu}
\end{figure}

Figure~\ref{fig:mu_xt} shows the space--time evolution of the optimized viscosity
$\mu(x,t)$. Similarly, figure~\ref{fig:mu_at_T} shows the normalized learned $\mu$ at T = 0.15 s.  The learned viscosity is highly localized near the moving discontinuities
associated with the hat-function edges, while remaining small in smooth regions.
Notably, the optimized viscosity exhibits \emph{sign-indefinite} behavior within the
allowed bounds, indicating that it does not act purely as a diffusive term. Instead,
$\mu(x,t)$ serves as a numerical correction that compensates for both dispersive and
dissipative errors of the central discretization. In the present experiment, the learned viscosity remains predominantly positive
($\mu_{\max}\approx9.5\times10^{-2}$), while weakly negative values
($\mu_{\min}\approx-5\times10^{-3}$) arise only locally near discontinuities,
providing controlled anti-diffusive corrections without compromising stability.

\subsection{Accuracy of the numerical solution}

\begin{figure}[H]
\centering

\begin{subfigure}{0.48\textwidth}
    \centering
    \includegraphics[height=4.5cm]{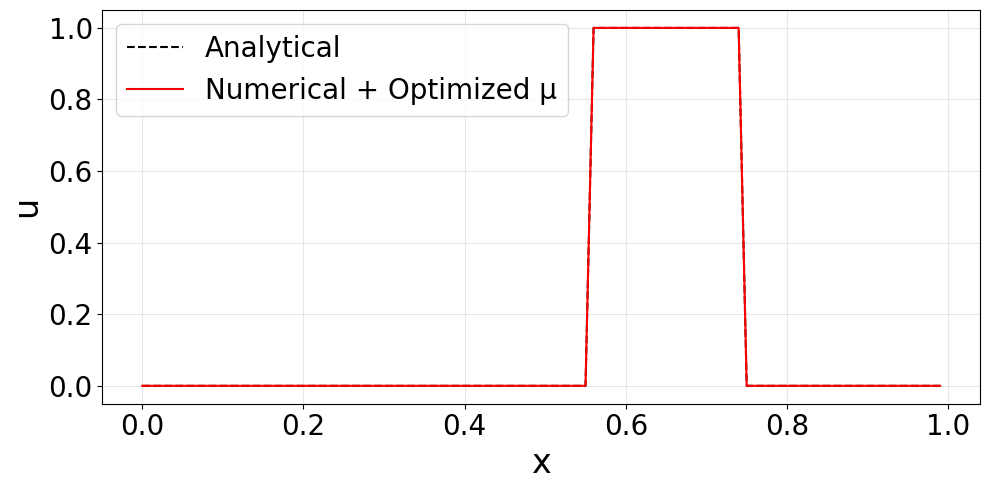}
    \caption{Solution at T = 0.15 s.}
    \label{fig:LCE_sol}
\end{subfigure}
\hfill
\begin{subfigure}{0.48\textwidth}
    \centering
    \includegraphics[height=4.5cm]{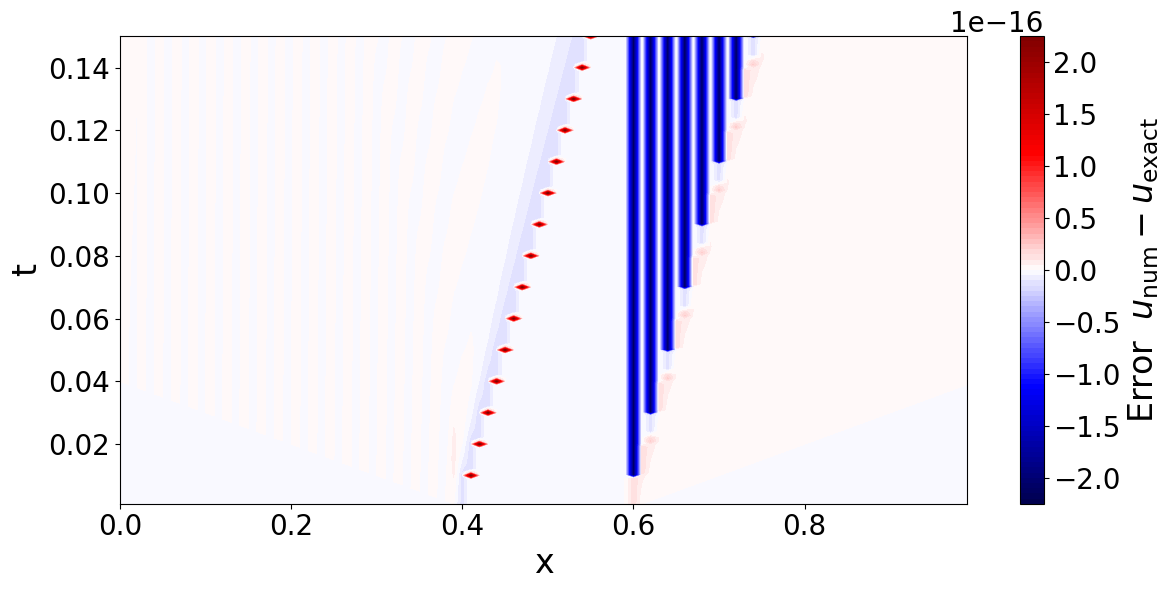}
    \caption{Error in the learned solution in $x-t$ plane}
    \label{fig:error}
\end{subfigure}

\caption{Comparison of learned artificial viscosity }
\label{fig:mu}
\end{figure}

Figure~\ref{fig:LCE_sol} compares the numerical solution obtained using the
optimized $\mu(x,t)$ with the exact solution at the final time. The error of the learned optimal viscosity solution in the $x-t$ plane is shown in figure~\ref{fig:error}.
Despite the known instability of the FTCS scheme for pure advection, the learned
viscosity stabilizes the method and enables accurate propagation of the discontinuous
profile.
Spurious oscillations typically associated with central schemes are strongly suppressed,
while excessive numerical smearing is avoided due to the localized nature of the learned
dissipation. The mean-squared error remains small throughout the simulation, confirming
the effectiveness of the data-driven viscosity closure.

\section{Discussion and Outlook}

The results presented in this work suggest a reinterpretation of artificial viscosity that unifies data-driven optimization, classical numerical analysis, and entropy stability within a single space--time closure framework.
Rather than acting as a purely dissipative spatial regularization, artificial viscosity emerges as a numerical reconvolution operator that compensates unresolved truncation errors introduced by discretization.
A central observation of this study is that gradient-based optimization, when applied to a fixed and consistent discretization, naturally identifies locally sign-indefinite viscosity fields.
The space--time closure perspective introduced here provides a mechanism-based explanation for this behavior and clarifies its relationship to established schemes such as Lax--Wendroff, Taylor--Galerkin, spectral viscosity, and compact flux corrections.

\section{Conclusion}

In this work, we investigated artificial viscosity from a data-driven,
space--time perspective using the linear convection equation as a model
problem.
By optimizing stabilizing closures for a fixed, consistent discretization,
we observed the systematic emergence of locally sign-indefinite viscosity
fields that nevertheless remain stable and accurate.
Although such behavior is not readily explained by classical
interpretations of artificial viscosity as a purely dissipative
mechanism, it motivates a reinterpretation of artificial viscosity as a
space--time closure that compensates for unresolved truncation errors
while preserving numerical stability.

Within this framework, the Lax--Wendroff scheme was shown to correspond to a degenerate time-to-space projection in which temporal truncation errors are eliminated and reintroduced as spatial diffusion.
Relaxing this restriction allows negative spatial viscosity to arise naturally as a numerical mechanism for compensating dispersive truncation errors, without violating global entropy stability.
The analysis further demonstrated that entropy stability constrains the integrated dissipation budget rather than the pointwise sign of individual operators.
This perspective provides a unifying conceptual framework for classical and data-driven stabilization methods and offers guidance for the design of adaptive, interpretable closures in future numerical schemes.
\section*{Future Scope}
Learning a consistent algebraic or simple closure for sign-indefinite artificial viscosity opens a pathway toward data-driven solvers that may substantially reduce the computational cost of DNS and LES for nonlinear conservation laws. Nevertheless, the discovery of stable, near-exact data-driven representations for inviscid systems such as the Euler and Burgers equations remains unresolved. Future efforts will explore whether physically constrained learning of artificial dissipation can bridge this gap and enable reliable modeling of inviscid dynamics.
\bibliographystyle{plain}
\bibliography{references}

@article{sengupta2023quantifying,
  title={Quantifying parameter ranges for high fidelity simulations for prescribed accuracy by Lax--Wendroff method},
  author={Sengupta, Tapan K and Suman, VK and Sengupta, Soumyo and Sundaram, Prasannabalaji},
  journal={Computers \& Fluids},
  volume={254},
  pages={105794},
  year={2023},
  publisher={Elsevier}
}

@article{NeelanDispersion2025,
  title={Improved approximate dispersion relation analysis using deep neural network},
  author={Neelan, Arun Govind},
  journal={International Journal of Computer Mathematics: Computer Systems Theory},
  volume={9},
  number={3},
  pages={155--182},
  year={2024},
  publisher={Taylor \& Francis}
}

@article{neelan2024physics,
  title={Physics-informed neural networks and higher-order high-resolution methods for resolving discontinuities and shocks: A comprehensive study},
  author={Neelan, Arun Govind and Krishna, G Sai and Paramanantham, Vinoth},
  journal={Journal of Computational Science},
  volume={83},
  pages={102466},
  year={2024},
  publisher={Elsevier}
}

@article{vonNeumann1950,
  title={A method for the numerical calculation of hydrodynamic shocks},
  author={VonNeumann, John and Richtmyer, Robert D},
  journal={Journal of Applied Physics},
  volume={21},
  number={3},
  pages={232--237},
  year={1950},
  publisher={AIP}
}

@article{Lax1960,
  title={Systems of conservation laws},
  author={Lax, Peter D and Wendroff, Burton},
  journal={Communications on Pure and Applied Mathematics},
  volume={13},
  number={2},
  pages={217--237},
  year={1960},
  publisher={Wiley}
}

@inproceedings{Jameson1981,
  title={Numerical solutions of the Euler equations by finite volume methods using Runge-Kutta time-stepping schemes},
  author={Jameson, Antony and Schmidt, Wolfgang and Turkel, Eli},
  booktitle={14th Fluid and Plasma Dynamics Conference},
  pages={1259},
  year={1981}
}

@article{Harten1983,
  title={High resolution schemes for hyperbolic conservation laws},
  author={Harten, Ami},
  journal={Journal of Computational Physics},
  volume={49},
  number={3},
  pages={357--393},
  year={1983},
  publisher={Elsevier}
}

@article{Tadmor1987,
  title={The numerical viscosity of entropy stable schemes for systems of conservation laws. I},
  author={Tadmor, Eitan},
  journal={Mathematics of Computation},
  volume={49},
  number={179},
  pages={91--103},
  year={1987}
}

@article{Donea1984,
  title={A Taylor-Galerkin method for convective transport problems},
  author={Donea, Jean},
  journal={International Journal for Numerical Methods in Engineering},
  volume={20},
  number={1},
  pages={101--119},
  year={1984},
  publisher={Wiley}
}

@article{Shu1988,
  title={Efficient implementation of essentially non-oscillatory shock-capturing schemes},
  author={Shu, Chi-Wang and Osher, Stanley},
  journal={Journal of Computational Physics},
  volume={77},
  number={2},
  pages={439--471},
  year={1988},
  publisher={Elsevier}
}

@article{Raissi2019,
  title={Physics-informed neural networks: A deep learning framework for solving forward and inverse problems involving nonlinear partial differential equations},
  author={Raissi, Maziar and Perdikaris, Paris and Karniadakis, George E},
  journal={Journal of Computational Physics},
  volume={378},
  pages={686--707},
  year={2019},
  publisher={Elsevier}
}

@article{BarSinai2019,
  title={Learning data-driven discretizations for partial differential equations},
  author={Bar-Sinai, Yohai and Hoyer, Stephan and Hickey, Jason and Brenner, Michael P},
  journal={Proceedings of the National Academy of Sciences},
  volume={116},
  number={31},
  pages={15344--15349},
  year={2019},
  publisher={NAS}
}

@article{Kochkov2021,
  title={Machine learning--accelerated computational fluid dynamics},
  author={Kochkov, Dmitrii and Smith, Jamie A and Alieva, Ayya and Wang, Qing and Brenner, Michael P and Hoyer, Stephan},
  journal={Proceedings of the National Academy of Sciences},
  volume={118},
  number={21},
  pages={e2101784118},
  year={2021},
  publisher={NAS}
}

@article{Duraisamy2019,
  title={Turbulence modeling in the age of data},
  author={Duraisamy, Karthik and Iaccarino, Gianluca and Xiao, Heng},
  journal={Annual Review of Fluid Mechanics},
  volume={51},
  pages={357--377},
  year={2019},
  publisher={Annual Reviews}
}

@article{Sirignano2018,
  title={DGM: A deep learning algorithm for solving partial differential equations},
  author={Sirignano, Justin and Spiliopoulos, Konstantinos},
  journal={Journal of Computational Physics},
  volume={375},
  pages={1339--1364},
  year={2018},
  publisher={Elsevier}
}

@inproceedings{Holl2020,
  title={Learning to control PDEs with differentiable physics},
  author={Holl, Philipp and Koltun, Vladlen and Thuerey, Nils},
  booktitle={International Conference on Learning Representations},
  year={2020}
}

@article{Baydin2018,
  title={Automatic differentiation in machine learning: a survey},
  author={Baydin, Atilim Gunes and Pearlmutter, Barak A and Radul, Alexey A and Siskind, Jeffrey Mark},
  journal={Journal of Machine Learning Research},
  volume={18},
  number={153},
  pages={1--43},
  year={2018}
}

@article{Warming1974,
  title={The modified equation approach to the stability and accuracy analysis of finite-difference methods},
  author={Warming, Robert F and Hyett, BJ},
  journal={Journal of Computational Physics},
  volume={14},
  number={2},
  pages={159--179},
  year={1974},
  publisher={Elsevier}
}

@article{Giles2000,
  title={An introduction to the adjoint approach to design},
  author={Giles, Michael B and Pierce, Niles A},
  journal={Flow, Turbulence and Combustion},
  volume={65},
  pages={393--415},
  year={2000},
  publisher={Springer}
}

@article{Jameson1988,
  title={Aerodynamic design via control theory},
  author={Jameson, Antony},
  journal={Journal of Scientific Computing},
  volume={3},
  pages={233--260},
  year={1988},
  publisher={Springer}
}

@article{Deconinck2013,
  title={Residual distribution schemes: foundation and analysis},
  author={Deconinck, Herman and Ricchiuto, Mario},
  journal={Encyclopedia of Computational Mechanics. John Wiley \& Sons, Ltd},
  year={2007}
}

@article{Tadmor1989,
  title={Convergence of spectral methods for nonlinear conservation laws},
  author={Tadmor, Eitan},
  journal={SIAM Journal on Numerical Analysis},
  volume={26},
  number={1},
  pages={30--44},
  year={1989},
  publisher={SIAM}
}

@article{Lele1992,
  title={Compact finite difference schemes with spectral-like resolution},
  author={Lele, Sanjiva K},
  journal={Journal of Computational Physics},
  volume={103},
  number={1},
  pages={16--42},
  year={1992},
  publisher={Elsevier}
}

@article{Guermond2011,
  title={Entropy viscosity method for nonlinear conservation laws},
  author={Guermond, Jean-Luc and Pasquetti, Richard and Popov, Bojan},
  journal={Journal of Computational Physics},
  volume={230},
  number={11},
  pages={4248--4267},
  year={2011},
  publisher={Elsevier}
}

@book{Hirsch2007,
  title={Numerical computation of internal and external flows: The fundamentals of computational fluid dynamics},
  author={Hirsch, Charles},
  year={2007},
  publisher={Elsevier}
}

@article{Beck2019,
  title={Machine learning approximation algorithms for high-dimensional fully nonlinear partial differential equations and second-order backward stochastic differential equations},
  author={Beck, Christian and E, Weinan and Jentzen, Arnulf},
  journal={Journal of Nonlinear Science},
  volume={29},
  number={4},
  pages={1563--1619},
  year={2019},
  publisher={Springer}
}

@article{Pathak2022,
  title={Fourcastnet: A global data-driven high-resolution weather model using adaptive fourier neural operators},
  author={Pathak, Jaideep and Subramanian, Shashank and Harrington, Peter and Raja, Sanjeev and Chattopadhyay, Ashesh and others},
  journal={arXiv preprint arXiv:2202.11214},
  year={2022}
}

@article{Tadmor2003,
  title={Entropy stability theory for difference approximations of nonlinear conservation laws and related time-dependent problems},
  author={Tadmor, Eitan},
  journal={Acta Numerica},
  volume={12},
  pages={451--512},
  year={2003},
  publisher={Cambridge University Press}
}

@article{Fjordholm2012,
  author  = {Fjordholm, Ulrik S. and Mishra, Siddhartha and Tadmor, Eitan},
  title   = {Arbitrarily High-Order Accurate Entropy Stable Essentially Nonoscillatory Schemes for Systems of Conservation Laws},
  journal = {SIAM Journal on Numerical Analysis},
  volume  = {50},
  number  = {2},
  pages   = {544--573},
  year    = {2012},
  doi     = {10.1137/110836961}
}
\end{document}